\documentclass[11pt]{article}
\usepackage{geometry}                
\usepackage{graphicx}
\usepackage{amssymb}
\usepackage{color}

\definecolor{mygold}{rgb}{0.9,0.9, 0.1}

\usepackage{epstopdf}

\usepackage{graphicx}
\usepackage{amssymb}
 \usepackage{psfrag}
 \usepackage{color}
\usepackage{amsthm}
\newtheorem{thm}{Theorem}
\newtheorem{cor}{Corollary}
\newtheorem{lem}{Lemma}
\newtheorem{prop}{Proposition}
\theoremstyle{definition}
\newtheorem*{defn}{Definition}
\theoremstyle{remark}

\def\proof{\noindent{\it Proof:}\ }
\usepackage{epstopdf}
\def\Bbb{\mathbb}
\def\R{{\Bbb R}}

\def\C{{\Bbb C}}
\def\Z{{\Bbb Z}}
\def\E{{\Bbb E}}
\def\H{{\Bbb H}}

\def\Q{{\Bbb Q}}
\def\U{{\Bbb U}}

\def\mbf0{\mathbf{0}}

\newcommand{\bmat}{\begin{pmatrix}}
\newcommand{\emat}{\end{pmatrix}}

\date{}

\title{Cusped  hyperbolic  3-manifolds: canonically CAT(0) \\ with  
  CAT(0) spines}

\author{Iain Aitchison}

\begin{document}
\maketitle

 \begin{abstract}
 
 We prove that every finite-volume hyperbolic 3-manifold $M^3$ with $p\geq 1$ cusps
 admits a canonical, complete, piecewise Euclidean CAT(0) metric, with a canonical projection to a CAT(0) spine
 $K_M^*$. 
Moreover:  
 \begin{description}

\item[(a)] The universal cover of $M^3$ endowed with the CAT(0) metric is a union of Euclidean half-spaces, glued together by identifying  Euclidean polygons in their bounding planes by pairwise isometry; 

\item[(b)] Each cusp of $M^3$ in the CAT(0) metric is a non-singular 
 metric product  $\E t_i^*\times [1,\infty)$, where 
 $ \{ \E t_i^*\}_{i=1}^p $ is a set  of  Euclidean cusp  tori,    with $\E t_i^*$ having the  canonical shape associated with the $i$th cusp;
 
\item[(c)] Metric singularities are concentrated on the 1-skeleton of $K_M^*$
with cone angle  $k\pi$ on any edge   of degree $k$.
The CAT(0) 2-complex $K_M^*$  is  constructed canonically 
from Euclidean polygons $ P^e_{i,j}$, 
which   reassemble to create  $ \{ \E t_i^*\}_{i=1}^p $; 

\item[(d)]  There is a canonical 1-parameter metric deformation, through piecewise-constant-curvature complete metrics, from the hyperbolic metric with limit the piecewise Euclidean one (facilitated by a simple application of
Pythagorus' Theorem);

\item[(e)] The hyperbolic metric on $M$ can be reconstructed from  
 a finite set of  points $p_{i,j}$
on the tori $\E t_i^*$, weighted by real numbers $w_{i,j} \in (0,1)$.

\end{description}

Our CAT(0) construction  can be considered `dual' to  
that of Epstein and Penner, but uses much simpler arguments, directly and canonically 
based on Ford domains.
Epstein and Penner's metrics, parametrized by a choice ${\mathcal T}$ of disjoint cusp horotori,
gives rise to  incomplete piecewise Euclidean metrics with singularities in cusps.
To each such choice ${\mathcal T}$, we also construct  a complete  CAT(0) metric of 
the above form, with CAT(0) spine $K_{{\mathcal T}}$. This CAT(0) metric structure is already visible via both Weeks' Snappea program, and its recent manifestation SnapPy by Culler and Dunfield, although its existence  has not previously been observed. 

Our  construction also   generalizes  to finite-volume $p$-cusped $n$-manifolds $W^n$, to  endow each  with a complete piecewise-Euclidean CAT(0) metric
with non-singular product end structures, whose singularities are concentrated in codimension 2:  such $W^n$ deformation retract to a natural spine, which is CAT(0) 
 as a manifestation of  polar duality of ideal hyperbolic polytopes.  

\end{abstract}

{\footnotesize
2010 MSC. 
Primary:   57M50,   20F67 ;  
Secondary:  51M10, 52B70  
}

\vfill
\pagebreak
\section{Introduction}\label{sec:intro}
	
	Spaces of constant curvature     play a  fundamental role in pure
mathematics, since their advent as solutions to
the conceptual problem of
the independence of Euclid's 5th Postulate. 
Relationships between 
geometry, complex analysis  and number theory, and between 
discrete and continuous
representations of mathematical objects, 
continue to be of profound
significance. We demonstrate
a 
simple interplay of the 
combinatorics and smooth structure
of non-compact hyperbolic space-forms,   geometrically dual 
to and inspired by Epstein and Penner's  canonical 
decompositions 
\cite{EP}. 
 
Generically, a closed hyperbolic manifold $M$ does not have a natural distinguished finite set 
of points by which to create a combinatorial structure: any given finite set of points enables the construction of a Delauney decomposition, or dually, a Dirichlet/Voronoi decomposition of $M$ into hyperbolic cells, as demonstrated by N\"a\"at\"anen-Penner
\cite{NP} in two dimensions, and more generally by
Charney-Davis-Moussong \cite{CDM},  using the Minkowski model for
hyperbolic geometry. 
 The results of \cite{CDM} also pertain to non-compact hyperbolic manifolds, but are 
less canonical in that case, requiring a careful choice of infinitely many points by which to create a locally-finite cell decomposition into compact cells.   

In these constructions, it can be shown that hyperbolic cells can be replaced by Euclidean ones: $M$ admits  
piecewise Euclidean CAT(0) structures, in the sense of Gromov. 
For dimension greater than 2, this utilises  
Rivin's description of ideal hyperbolic polyhedra in 3-space,
using  polar duality, and its generalisation to higher dimensions  by  Charney-Davis \cite{CD}.

The results of \cite{CDM, NP} were motivated by Epstein--Penner's \cite{EP}
canonical  decomposition of a cusped hyperbolic manifold $M$ of finite volume
into ideal hyperbolic polytopes, using their 
convex hull construction in Minkowski space: this  construction exploits
the distinguished finite set of `ideal points' of $M$ corresponding to   cusps. 
The Epstein--Penner piecewise Euclidean metrics arise naturally  by replacing hyperbolic polytopes with Euclidean ones, but  are incomplete, and have singular set 
intersecting cusps. 

Epstein--Penner's canonical decomposition    is essentially
a  Delauney decomposition, based on the  ideal points of $M$. Dual to any Delauney decomposition is a Voronoi--Dirichlet  decomposition, and Epstein and Penner show that their decomposition of $M$
into hyperbolic polytopes is naturally dual to the classical Ford decomposition, which is traditionally defined using isometric circles in the  upper-half-space model for hyperbolic space.  Our proof is based on the geometry of Ford domains, viewed -- as heuristically described in \cite{EP} -- as arising by the collision of expanding horospheres, and is thus based on a Voronoi--Dirichlet  construction.

Thus, of all four CAT(0) structures defined on a cusped hyperbolic manifold $M$, arising respectively from the Delauney, Voronoi--Dirichlet, Epstein--Penner, and Ford decompositions of $M$  into hyperbolic pieces, the latter two are most 
natural, but only our decomposition gives a complete metric, a CAT(0) spine, and singularity-free cusps.

The proof of our main theorem generalises to any dimension,
and so we concentrate on dimension $3$ for the purposes of illustration and exposition, and the contextual significance of the conclusion: for 2-dimensional analogues, with applications to Riemann surfaces and moduli thereof, we refer to the work of Bowditch and Bowditch --Epstein \cite{Bo, BE}. For higher dimensions, the proof is essentially identical, given the
Charney--Davis  results on polar duality \cite{CD}.

\section {Acknowledgement}\label{sec:ackno}

The author thanks  Igor Rivin, Norman Wildberger, Dani Wise, Brian Bowditch,  Makoto Sakuma, Jeff  Weeks and Mladen Bestvina
for inspiration, valuable conversations, suggestions,  and encouragement.
Analogous results for alternating links were mentioned at the 
Osaka `Knots '90'
conference, and were described in detail, in the context of cubings
of manifolds, in several
lectures in Tokyo,  Osaka and Kobe in 1992, and the author 
would like to express his thanks for 
his experience of Japanese hospitality.
The observation of the existence of CAT(0)  cubical metrics owes some 
debt to a question asked of the
author by Darren Long in 1990, concerning the $8_{17}$ knot. 
That  CAT(0) metrics are so
ubiquitous, for all cusped hyperbolic 3-manifolds was realised 
subsequent to a very helpful  conversation with Dani  Wise in 1999. 

\section{Preliminaries and notation}\label{sec:prelim}

We denote hyperbolic 3-space by $\H^3$, and by $\bar \H^3$ its compactification obtained by adding ideal points.
These constitute the sphere at infinity, $S^2_\infty = \bar \H^3 -\H^3$.
The upper-half space model  $\U\H^3$ for $\H^3$  has underlying set  
$$
\U\H^3 := \{ (x,y,z)\in \R^3 \, | \, z>0\} = \R^3_+ $$
and  sphere at infinity represented as 
$$S^2_\infty = \R^2_0 \cup  \infty  :=   \{ (x,y,0) \in \R^3 \}  \cup  \infty,
$$ 
with $\R^2_0$ inheriting a Euclidean metric, up to similarity: for any $p\in \R^2_0$, any dilation of $\R^3$
centered at $p$ gives a hyperbolic isometry, as does any translation fixing $\R^2_0$ setwise.

 In $\U\H^3$, any horosphere appears as  either a horizontal Euclidean plane,
 $$
{\it  HoP}_a := \{     (x,y, a) \in \R^3_+ \} ,$$
or as a Euclidean sphere $
{\it  HoS}_{p,d}$ of diameter $d$,   tangent to $\R^2_0$ at $p$, which is deleted. 

Similarly, in $\U\H^3$,  a hyperbolic plane appears as either a  vertical Euclidean plane,
$$
{\H yP^2}_{a,b,c} := \{     (x,y, z) \in \R^3_+ \, | \, ax+by+c=0 \}  \cup \infty,$$
after deleting  its  {\it circle-at-infinity} 
$
{S^1}_{a,b,c} := \{     (x,y, 0) \in \R^2_0 \, | \, ax+by+c=0 \}  \cup \infty;$
or as a Euclidean hemisphere $ {\H yS^2}_{p,r}$ of radius $r$, centered on $\R^2_0$ at $p$, with  the equatorial boundary circle $ {S^1}_{p,r}\subset\R^2_0$,  its {circle-at-infinity},   deleted.

If $P^h\subset  {\H yS^2}_{p,r}$ is any compact,   \emph{hyperbolic},  polygon, its orthogonal projection to $\R^2_0$
is a compact,   Euclidean polygon $P^e$, with respect to the standard Euclidean metric on $\R^2_0$. Each edge of $P^e$ determines a vertical plane in $\R^3_+$, and hence a hyperbolic plane in $\U\H^3$ intersecting $ {\H^2}_{p,r}$ in a hyperbolic geodesic containing a geodesic boundary segment of $P^h$.

\defn
A \emph{label} for $P$ is a pair $(p,r) \in \R^2_0\times \R_+\cong \C\times \R_+$.
(We do not associate a label to polygons lying in vertical hyperbolic planes.)

 A non-compact complete hyperbolic 3-manifold $M^3$ of finite volume 
$M^3$ has $p$ cusps for some $p\geq 1$, and decomposes
 \cite{Th1,Th2}    as 
$M^3 =  M^{thick}\  \cup\  \{ {\mathcal C}_i\}_{i=1}^p,  $
with compact  `thick' part $M^{thick}$ having as complement
a disjoint union of $p$  cusps ${\mathcal C}_i$, $i=1,\dots , p$, each topologically a product $T^2\times (1,\infty)$ of tori.
In the
hyperbolic metric, each  torus $T^2\times \{ t \}$ inherits  a Euclidean metric, whose scale shrinks  exponentially as 
$t \to \infty$ 
 at unit speed. Accordingly, there is canonically associated  to $M^3$ 
a set $\{ \E c_i\}_{i=1}^p$ of elliptic curves, with $  \E c_i $ associated  to the $i$th cusp ${\mathcal C}_i$.

\medskip

A \emph{cusp curve} for $M^3$  is any elliptic curve $ \E c_i$ associated to a cusp of $M^3$. A \emph{cusp torus}, or \emph{horotorus} in $M^3$ is the image of any \emph{Euclidean torus}   $ \E t$ isometrically embedded
in some cusp of $M^3$.

\medskip

Given any set  $\{ \E c_i\}_{i=1}^p$  of  elliptic curves, since each curve $\E c_i$ admits a unique flat Euclidean metric, up to scale, we obtain a set of Euclidean tori  $\{ \E t_i\}_{i=1}^p$ by independently specifying a scale for each.    
A priori, there is no
specified 
scale for each cusp elliptic curve: specifying a scale amounts to choosing a cusp torus, 
and each cusp determines a least upper bound for the possible size of any  cusp torus, and hence determines a 
\emph{maximal cusp torus} $\E^* t_i$, $i=1,\dots, p$. Each $\E^* t_i$ has a non-empty finite set of self-tangencies, and thus determines a finite set of points on the corresponding elliptic curve.

 \medskip

By Marden and Prasad's generalization \cite{Ma,Pr} of Mostow rigidity, there is a unique (up to
conjugation) representation $\rho : \pi_1 (M^3) \longrightarrow 
\Gamma = \rho (\pi_1 (M^3)) \subset PSL_2(C)\cong Isom_+(\H^3)$, with $M \cong \H^3/\Gamma$.
Thus $\Gamma$ naturally acts on $\U\H^3$, by M\"obius transformations on the Riemann sphere $\hat S^2 = \C\cup\infty \cong \R_0^2\cup \infty$, and 
$\Gamma$   acts on  
$\bar H^3 $ with a dense set ${\mathcal P}_\Gamma\subset S^2_\infty$ of 
parabolic fixed points,   falling into $p$
distinct orbits corresponding to the $p$ cusps of $M^3$.

\medskip

The preimage in $\H^3$ of any cusp torus is a disjoint set of
 horospheres   with  inherited  Euclidean metric. Such  metrics can be seen 
algebraically in the Lorentzian model, and are visually natural in the upper-half space model for horospheres 
${\it Hop}_a$ centered at $\infty\in  S^2_\infty = \infty \cup \R^2_0$. 
Accordingly the upper-half-space models we use to describe $M^3$ will have $\infty$ as a parabolic fixed point for $\Gamma$:
For $p_i\in S^2_\infty$ a parabolic fixed point corresponding to the $i$th cusp, we conjugate $\Gamma$ to 
$\Gamma_i$ so that
$p_i$ is at $\infty_i := \infty$ in $\U\H^3_i := \U\H^3$; thus the stabilizer of $\infty_i$ is a $\Z\oplus \Z$ subgroup of $\Gamma_i$.
When $p>1$,   although  $\Gamma_i$ is   conjugate to $\Gamma_j$ in $PSL(2,\C)$, 
there is no  element of $\Gamma_i$ sending $\infty $ to any parabolic fixed point corresponding to a distinct cusp 
${\mathcal C}_j$ of $M^3$, 
since the images of $\infty$ under $\Gamma_i$  constitute  a single orbit of parabolic fixed points.

\defn Denote by  $\U\H^3_i$ the $i$th upper-half-space on which $\Gamma_i$ acts. When a horotorus $\E t_i \subset {\mathcal C}_i$ is specified, 
  we will generally assume that we have conjugated $\Gamma$ so that  the horosphere $HoP_1$ projects  to  $\E t_i  $ under the action of $\Gamma_i$ on $\U\H^3_i$. The metric on the   cusp torus $\E t_i $ is determined by the action of the parabolic subgroup stabilizing $\infty$, since the induced Euclidean metric on $HoP_1=\{ (x,y,1)\} \cong \{ (x,y)\} =  \R^2$ is the standard Euclidean metric.

\section{Statement of results}\label{sec:res}

\begin{thm}
  Suppose $M^3$ is a non-compact,  
connected  
$3$-manifold admitting a complete hyperbolic 
metric of finite-volume with $p\geq 1$-cusps. Then
\begin{itemize}

\item $M^3$ admits 
a  complete piecewise-Euclidean CAT(0)  metric, 
with singular set concentrated on a finite connected graph; all edge cone angles are  of form 
$k\pi , \ 3\leq k\in \Z$.

\item Each hyperbolic  cusp ${\mathcal C}_i, \ i= 1,\dots , p$, 
canonically 
determines an elliptic curve  $\E c_i$ (a  Euclidean similarity class of a closed 
Euclidean torus): 
there is a consistent 
choice $\E t^*_i$
of a representative Euclidean torus from each class, inducing the
CAT(0) metric on $M^3$, with each cusp   the Euclidean metric  product 
$ \E t^*_i \times [1,\infty )$.

\item Each 
$\E t^*_i$ is a
union of    convex
Euclidean  polygons ${}P^e_{i,j}$:   the  
CAT(0) metric  on $M^3$
arises as the quotient space of 
$\coprod  \E t^*_i \times [1,\infty )$
by Euclidean isometric identification 
of pairs of polytopes ${} P^e_{i,j}
\times
\{ 1
\}$.

\item The piecewise-Euclidean 2-complex  $K$ obtained 
by
pairwise identification of polytopes 
$ P^e_{i,j} \times \{ 1 \}$ 
 is a spine for $M^3$, and is CAT(0). 
 
 \item There is   canonical deformation between 
the unique hyperbolic and CAT(0)  metrics, via a natural manifestation of Pythagorus' Theorem.

 \item The decomposition of Euclidean tori into polygons $ P^e_{i,j}$  is determined by a canonical  finite set  of weighted points $p_{i,j}
\in \E t^*_i\times (0,1)$. The hyperbolic metric on $M^3$ can be reconstructed from the data 
$\{p_{i,j}\}$.

\item The universal cover of $M^3$ with CAT(0) metric is a  union of Euclidean half-spaces, corresponding to  hyperbolic horoballs,
glued together by pairwise isometry of Euclidean polygons forming tessellations of their bounding Euclidean planes.  

\end{itemize}
\end{thm}

\section{Ford domains}\label{sec:ford}
 Epstein and Penner  define their  canonical decomposition using the Lorentzian model: traditional Ford domains are   naturally seen in the upper-half space model, since
  classically  their construction is via isometric circles in $\R_0^2$. 
In the 1-cusped case, Epstein and Penner formalise the 
heuristic `bumping locus' construction  by expanding horospheres in $\U\H^3$ \cite{EP}, showing that their canonical decomposition of $M^3$ into ideal polyhedra   is naturally dual to the Ford complex:
we describe  `horospherical bumping' for $p\geq 1$ in more detail, working equivariantly in 
$\U\H^3$ as universal cover of $M^3$:

Take any disjoint union ${\mathcal T} = \{ \E t_i\}_{i=1}^p\subset M^3$ of cusp tori, one for each cusp.
These lift to a union ${\mathcal H} = \{ H_p\}$ of horospheres  
centered 
at parabolic
fixed points $p\in S^2_\infty$, equivariant with respect to the action of $\Gamma$, and  
determine a set of disjoint open horoballs  
 ${\mathcal B} = \{ B_p\}$.
Expand each $H_p$ at unit speed, allowing them to `flatten' against each other, creating a locally finite piecewise geodesic 2-complex $K_{\mathcal H}$. 

\defn We call this  `bumping locus'  $  K_{\mathcal H}$ the 
\emph{Ford complex} for $\Gamma$ determined by   ${\mathcal H}$. Its projection $K_{\mathcal T} =  K_{\mathcal H}/\Gamma$ is a 2-complex in $M^3$, which is the \emph{Ford spine} for $M^3$ determined by ${\mathcal T}$. There is a strong deformation retraction from $M^3$ to $K_{\mathcal T}$.

 \medskip

  Viewed from  $\infty_i$  for $\U\H^3_i$, the visible part $  K^{\infty_i}_{\mathcal H}$ of $  K_{\mathcal H}$ is a locally finite piecewise geodesic 2-complex
 constructed from compact hyperbolic polygons $P^h_{i,j}$,   projecting to Euclidean polygons $P^e_{i,j}$  tessellating 
$HoP_1$  (or, equivalently, $\R^2_0$).   
 
 \medskip

 The Ford complex is created by numerous expanding geodesic discs in such hyperplanes, which in turn intersect each other creating the 1-skeleton: viewed from $\infty_i$, the Euclidean projections of these discs expand until they encounter  other expanding discs, at which stage their boundary circles also `flatten' against each other creating the straight boundary-edges of Euclidean polygons $P^e_{i,j}$. However, these expanding Euclidean discs do not expand at constant rate: we   discuss this later.

 \medskip

 \defn By abuse of language, we call the closure  of the complementary component of  $  K^{\infty_i}_{\mathcal H}$ in $\U\H^3_i$ containing the horoball $B_\infty$ a \emph{Ford ball} for $\Gamma_i$, denoted    by   $FB_{{\mathcal H},i}$: this is non-compact,   has boundary $  K^{\infty_i}_{\mathcal H}$ with infinitely many faces, and is the analogue of a Dirichlet polyhedron, with center at $\infty_i$.
 
 \medskip
 The pair ($\Gamma,  {\mathcal H})$ determines an equivariant tessellation  of $\H^3$ by copies of Ford balls $FB_{{\mathcal H},i}$.  This tessellation is equivalently created  by uniformly expanding  all horoballs in ${\mathcal B}$,
 allowing them to flatten against each other. Each $FB_{{\mathcal H},i}$ is stabilized by a $\Z\oplus\Z$ subgroup, and projects to a neighbourhood, denoted   $FC_{{\mathcal K},i}$, of the cusp ${\mathcal C}_i$. We call these
 \emph{Ford cusps}: these are the closures 
 in $M^3$ of the complementary components to $K_{\mathcal T}$, and each contains a horotorus $\E t_i$ naturally decomposed as a union of Euclidean polygons $P^e_{i,j}$. 
 
 \begin{prop}
 $\H^3$ is obtained from the disjoint  union of  Ford balls,    by   isometric pairwise identification of hyperbolic polygons
 $P^h_{i,j}$ in their boundaries.
 \end{prop}

  \medskip
  
   \defn  Define a \emph{hyperbolic Ford polytope} $FP^h_{{\mathcal H},i,j}$ for $\Gamma_i$ to be the closure in $\U\H^3_i$ of 
   the region vertically above a hyperbolic  polygon $P^h_{i,j}$: similarly, define a 
   \emph{Euclidean Ford polytope} $FP^e_{{\mathcal H},i,j}$  to be the closure in $\R^3_i$ (underlying $\U\H^3_i$) of 
   the region vertically above a Euclidean  polygon $P^e_{i,j}\subset HoP_1$. Both can be construed as `cones of polygons to infinity'.
Each $FB_{{\mathcal H},i}$ is a union of hyperbolic Ford polytopes, equivariantly with respect to the action of $\Gamma_i$. Thus

 \begin{prop}
 $M^3$ is obtained from the disjoint  union of hyperbolic Ford polytopes,    by pairwise   isometric identification of hyperbolic polygons 
  in their boundaries.
 \end{prop}
 
The finitely-many boundary faces of a Ford polytope consist of some hyperbolic $m$-gon $P^h_{i,j}$, together with $m$ non-compact hyperbolic triangles with exactly one ideal vertex.

 \medskip
Dual to the Ford complex is a decomposition of $\U\H^3$ into ideal polyhedra, which generically are 
simplices: the ideal cell dual to a given 0-cell $x$ of the Ford complex is the convex hull
 of the set of parabolic points determined by the closest equidistant horospheres. These polytopes form the Epstein--Penner \emph{canonical decomposition} \cite{EP}, which is unique when $p=1$, but otherwise admits a parameter space of real dimension $(p-1)$ corresponding to the $p$ choices of disjoint horotori, up to simultaneous rescaling.
 Akiyoshi has  shown in \cite{Aki}  that a finite volume
hyperbolic manifold with multiple cusps admits finitely many combinatorial types of canonical cell
decompositions.

\section{The distinguished  Ford complex}\label{sec:canon}

The `dynamic' view  allows us to 
generalise, to an arbitrary $p$-cusped manifold,  the  heuristic described and made more  precise  in the 1-cusped case  by Epstein--Penner  \cite{EP}. When $p=1$, all choices of embedded horotorus
are equal after contraction or expansion, and the Ford complex forms from the instant any expanding horotorus first contacts itself, and thus no choice in the construction is possible: it is natural for the unique Ford complex  to be seen as arising from expanding balloons flattening against each other, or via isometric spheres.  
An appropriate  generalisation of a distinguished choice for a Ford complex 
when $M^3$ has $p > 1$   cusps is not immediately clear from this `balloon flattening' perspective, since 
there is a $(p-1)$-dimensional parameter space for spines arising from possible initial choices of disjoint cusp tori. 
Consider three distinct heuristic scenarios as expanding horotori encounter each other, with a view to adjusting an initial family of embedded horotori to  create a more natural one:

\begin{description}
\item[Flattening:] This is described above: locally, expanding horospheres flatten against each other. The expansion process stops when each point of each horosphere has encountered another;
\item[Immersed transition:] Instead of flattening against each other, allow horospheres to continue expanding, becoming immersed. Collectively the expanding horospheres eventually pass through all  points in the complement of their horoball, and create quite complicated intersection patterns;
\item[Domination and submission:] Partition the set of cusps into two subsets, and declare one subset to be dominant, the other submissive. When two dominant horospheres meet, they flatten against each other; 
submissive horospheres are pushed back into their cusps by expanding dominant horospheres. In $\U\H^3_i$, for $i$ submissive, a horizontal horosphere eventually rises, supported by tangency with expanding dominant horospheres. Freeze the evolution at some instant, reverse the process by shrinking all immersed horotori back to disjointly embedded ones, and then allow all to expand again,
but now all flattening against each other.

\end{description}

\begin{prop} There is a distinguished  1-parameter family $H^*$ of disjointly embedded cusp tori giving a corresponding unique Ford complex $K^*_M$  for $M^3$, defined independently of  $p\geq 1$, naturally generalising the case $p=1$, by 
 viewing  the Ford complex as created by initially  intersecting expanding horospheres.
\end{prop}

\noindent{\bf Proof}: 
Consider expanding one arbitrarily  chosen embedded horotorus $\E t_i$, ignoring all others, and allowing self-intersection rather than self-flattening.  After finite time, the expanding torus $\E t^t_i$ sweeps past all
points of $M^{thick}$, leaving only points of cusps not yet encountered: the torus $\E t^t_i$ must, by this stage, be immersed, not embedded: this is clear when $p=1$; and for  $p=1$, clear since  the set of unencountered points of $M^3$ is disconnected.  

In keeping with the democratic philosophy of treating all cusps equally, we consider the set
$\{ \E t^t_i \}$ of all maximal cusp tori. Collectively they form a non-transverse immersion of 
$p$ Euclidean tori, and each of these can be independently shrunk at unit speed to become embedded. Doing this uniformly 
and simultaneously for each, we obtain a  regular homotopy of $p$ immersed Euclidean tori,
which eventually becomes  a  family ${\mathcal T}^*$ of embedded cusp tori
(with  1-parameter set of   choices by re-scaling uniformly). 
  Allowing all  these now-embedded cusp tori to expand again, we construct their collision locus $K_M^* := K_{{\mathcal T}^*}$.

\defn The \emph{distinguished} Ford spine {$K^*_M$ } for $M^3$ is defined to be
$K^*_M := K_{{\mathcal T}^*} $, for any such embedded family of rescaled maximal tori.

\medskip

 The defining characteristic for 
 {$K^*_M$ }  is arguably the most natural definition for distinguishing a family of cusp tori: others, perhaps less natural, can be defined using similar notions. For example, let $\E t^{thick}_i$ denote the immersed torus obtained by 
 expanding  $\E t^*_i$ until the first instant it has encountered each point of each
 other $\E t^*_j$. By this time, all points of $M^{thick}$ have been encountered by the $i$th expanding cusp torus.  Take the union $\{ \E t^{thick}_i \}_{i=1}^p$ of all these immersed tori, and uniformly and simultaneously shrink each backwards until each is embedded, and then allow the resulting embedded family ${\mathcal T}_{thick} $ to expand to create 
 $K^{thick}_M := K_{{\mathcal T}_{thick}}$, which is another natural choice for a distinguished family when $p\geq 1$.

\section{ 
Piecewise Euclidean structures: existence}\label{sec:euclid}

We now  define the Euclidean structure on $M^3$, arising from any Ford spine $K_{\mathcal T}$:

\defn The piecewise Euclidean structure $M^3_{\mathcal T}$ corresponding to
${\mathcal T}$ is defined by  
replacing each polygon $P^h_{i,j}$ by its projected Euclidean polygon $P^e_{i,j}$, 
 replacing each Ford polytope $FP^h_{{\mathcal H},i,j}$ by  the Euclidean Ford polytope
$FP^e_{{\mathcal H},i,j} := P^e_{i,j} \times [1,\infty)$, and replacing 
each Ford ball $FB_{{\mathcal H},i}$   by the Euclidean half-space
$FB^e_{{\mathcal H},i} \cong \R^3_{\geq 1} := \{ (x,y,z)\, | \  \, z\geq 1\}$, which is the union of Euclidean Ford polytopes.

\medskip

Heuristically  we vertically project those  hyperbolic polygons in 
$  K^{\infty_i}_{\mathcal H}$,  whose interior is  visible from $\infty_i$, to 
the horizontal plane at height 1, and take the vertical half-infinite prism above their images.
This is essentially coning each such hyperbolic polygon to the parabolic fixed point at
infinity,  intersecting with the horosphere at height 1, and `opening up' each Ford polytope
by no-longer-exponentially-shrinking Euclidean cusp geometry as we approach $\infty$. The labels assigned to polygons define an orthocenter for  both  
hyperbolic   and corresponding Euclidean polygons, as viewed from $\infty$.

It is clear that each Euclidean Ford ball $FB^e_{{\mathcal H},i}$ has a geometric combinatorial structure equivariant under the $\Z\oplus\Z$ action stabilizing $\infty_i$: we must show that the corresponding nonsingular quotients, which replace Ford cusps by Euclidean products,   glue together to produce the singular piecewise Euclidean structure $M^3_{\mathcal T}$, whose universal cover is then a geometrically complete union of Euclidean half-spaces..

\begin{thm}
Suppose $M^3$ is any non-compact,  
connected  
$3$-manifold admitting a complete hyperbolic 
metric of finite-volume with $p\geq 1$-cusps, with any specified  
complete family of disjoint horotori ${\mathcal T} =\{ \E t_i \}$.

\begin{itemize}

\item  The metric structure $M^3_{\mathcal T}$ is 
a  complete piecewise-Euclidean  metric, 
with singular set concentrated on a finite connected graph; all edge cone angles are  of form 
$k\pi , \ 3\leq k\in \Z$.

\item Each 
$\E t_i$ is decomposed naturally as a
union of  
Euclidean  polygons ${}P^e_{i,j}$:   the  
piecewise Euclidean  metric  $M^3_{\mathcal T}$
arises as the quotient space of 
$\coprod  \E t_i \times [1,\infty )$
by Euclidean isometric identification 
of pairs of polytopes ${} P^e_{i,j}
\times
\{ 1
\}$.

\item The piecewise-Euclidean 2-complex  $K_{\mathcal T}$ obtained 
by
pairwise identification of polytopes 
$ P^e_{i,j} \times \{ 1 \}$ 
 is a piecewise Euclidean spine for $M^3_{\mathcal T}$.

\end{itemize}
\end{thm}

\proof
Each hyperbolic polygon ${}P^h_{i,j}\in K^{\infty_i}_{\mathcal H}$   is contained in a unique hyperplane 
$HyS_{q, d}$: we assign the  label $(q,d) \in  \R^2_0\times \R_+\cong \R^3_+$   to  ${}P^h_{i,j}$. Thus $q$ is a parabolic fixed point for $\Gamma_i$, and ${}P^h_{i,j}$ lies in the hyperbolic plane formed by $H_{\infty_i}$  and 
$\H_q$  flattening against each other. Let $d = e^{-t}$, where $t$  denotes the time of initial tangency between these expanding horospheres since expansion began. 
Now $q$  lies in some orbit corresponding to a cusp ${\mathcal C}_k$, $k\in \{ 1, \dots , p\}$, and we consider the corresponding model $\U\H^3_k$. In this picture, some   polygon ${}P^{h}_{k,s}$ in the orbit of  ${}P^h_{i,j}$, and hence isometric to it by an orientation-reversing isometry (cf inversion in isometric spheres),
 is visible from $\infty_k$. These two hyperbolic polygons project to polygons in the boundary of Ford cusps in $M$, and are identified there by hyperbolic isometry gluing part of the boundaries of these cusps together.

\begin{figure}[htbp] 
   \centering
   \includegraphics[width=4.8in]{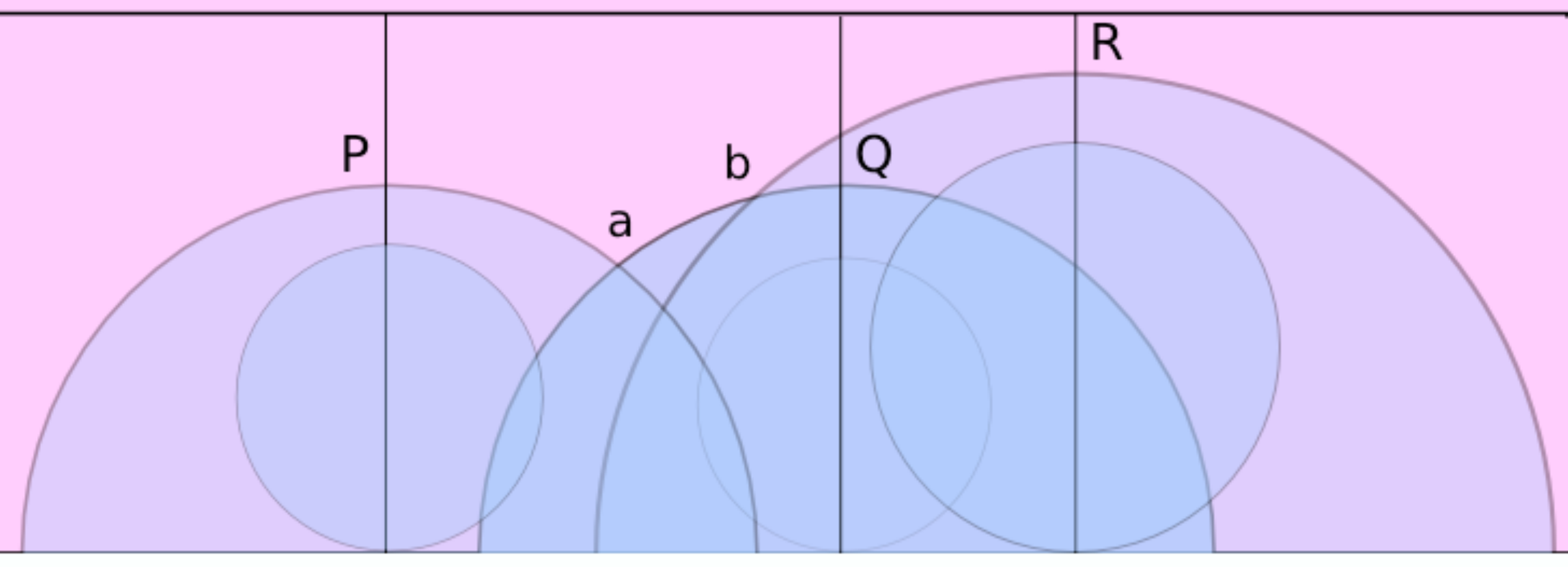} 
   \caption{Three expanding horospheres $HoS_*$, and one descending $HoP_*$: three hyperplanes
   $HyS_*$ with maximal points $P,Q,R$. The hyperbolic polygon $ab$ is part of the Ford complex: its label is 
   determined by the coordinates of $Q$ in $\R^3$. Note that $Q$ is invisible from $\infty$: $ab$ is contained in a hyperbolic disc centered at $Q$, but does not itself contain $Q$. }
   \label{fig: horolap}
\end{figure}

Consider the label $(q',d')$ for  ${}P^{h}_{k,s}$. Since $d' = e^{-t'}$ records the time of first tangency of expanding horospheres,
$t'=t$ and so $d'=d$: the corresponding hyperplanes in $\U\H^3_i, \U\H^3_k$ appear with the same Euclidean diameter. Similarly, the 
points of ${}P^{h}_{k,s}$ and ${}P^h_{i,j}$ are created by circle expansion, and so can be put in correspondence:
we may place both hyperplanes and polygons in the same $\U\H^3$, and observe they can be made to coincide by 
orientation-reversing Euclidean congruence of $\R^2_0$.

Summarizing:  if two hyperbolic polygons are identified by  an
element of $\Gamma$, the polygons each have the same height in their
corresponding rescaled half space models. But a hyperbolic polygon with given label
$(q,*)$
 in the upper half space model uniquely determines a Euclidean
similarity class of Euclidean polygons by vertical projection; and specifying the height $*$ uniquely
determines the scale.

Now take a disjoint union $\bigcup  FB^e_{{\mathcal H},i}$ of a countably infinite 
number of copies of each Euclidean Ford ball/half-space. Then  for each $i$, 
$\partial \, FB^e_{{\mathcal H},i}$ is a union of Euclidean polygons $P^e_{i,j}$, and to each we isometrically identify a corresponding $FB^e_{{\mathcal H},k}$ by isometric identification with $P^e_{k,s}$. The resulting 3-complex is homeomorphic to $\U\H^3$, and is metrically complete. 

The Ford complex $K_{\mathcal H}$ is 
replaced by, and is combinatorially equivalent to,  its piecewise Euclidean counterpart $K^e_{\mathcal H}$, 
obtained from the disjoint union 
 $\bigcup  \partial FB^e_{{\mathcal H},i}$ by pairwise isometric identification of all such Euclidean polygons $P^e_{i,j}$, $P^e_{k,s}$. Edges of  $K^e_{\mathcal H}$ correspond 
to those of $K_{\mathcal H}$, which have degree $\geq 3$ (generically each edge has exactly 3 polygons incident with it). Since the Euclidean edges lie in the boundaries of half spaces, all edges have  cone angle a multiple of $\pi$ in  $M^3_{\mathcal T}$.

Combinatorially, the structure is identical to that of the Ford complex, and so is  equivariant with respect to  the natural action of $\pi_1(M^3)$. The metric structure is equivariant, and so  descends to define the metric $M^3_{\mathcal T}$ with properties as stated in the theorem.

  The structure we describe is  no longer compatible with
representations of $\pi_1(M)$  in $PSL(2,\C)$, for all $M^3$, simultaneously acting as   isometries of the same space $\H^3$:   
Each piecewise Euclidean structure on $M^3$ endows its universal cover, topologically $\R^3$, with piecewise Euclidean  metrics which generally  differ for different $M^3$, and different choices for ${\mathcal T}$.

\section{Piecewise Euclidean structures are  CAT(0)}\label{sec:excat}

For basic definitions   for this section, we refer to \cite{BH, Ri1,Ri2, CD, CDM}. 
In order to prove that the piecewise Euclidean structures we have defined on   
$M^3_{\mathcal T}$ and  its spine $K^e_{\mathcal T}$ are CAT(0), we must argue that the link of each point is CAT(1): all  geodesics in each piecewise spherical link should be of length at least $2\pi$.
Such a piecewise spherical link is called \emph{large}:
there is a unique geodesic between any two points of distance less than $\pi$.
The essence of the argument is that the Ford complex is geometrically dual to the canonical Epstein--Penner canonical decomposition into finite-volume ideal hyperbolic polyhedra, and that the links of vertices in the piecewise Euclidean structures $K^e_{\mathcal T}$  and $M^3_{\mathcal T}$ are essentially the polar duals of these hyperbolic polyhedra. 
Rivin \cite{Ri2}   showed that the polar dual of a convex ideal hyperbolic  polyhedron is large in dimension 3. This result was generalized by Charney and Davis \cite{CD} to higher dimensions, and accordingly we adapt some of their notation so that relevant parts of their description are clearer in the present context. 

We must consider  the link of any point $x\in M^3_{\mathcal T}$ in the interior of a $k$-cell of $M^3_{\mathcal T}$, $k= 0,1,2,3$. Heuristically the metric for $M^3_{\mathcal T}$ should be CAT(0), since it is already so for the hyperbolic metric on $M^3$, where all links are then standard 2-spheres.   The solid angles at 0-cells created by intersecting with $K^h_{\mathcal T}$ are enlarged in $M^3_{\mathcal T}$,   becoming hemispheres: this should not create shorter geodesics. We describe the link structure with a little more care, since the more delicate structure of links in  $K^e_{\mathcal T}$ is also revealed.
It is important to note that the metric  2-complex $K^h_{\mathcal T}$ is not CAT(0), since the links of 0-cells are not large. 

For $k =3$, $x$ is an interior point of a Euclidean half space, and its link is thus a standard round sphere, with all geodesics of length $2\pi$.  For $k = 2,$ $x$ lies in the interior of some Euclidean polygon $P^e_{i,j}$, and its link in $M^3_{\mathcal T}$ is a union of two hemispheres corresponding to the two 
half-spaces identified along $P^e_{i,j}$, and again is a standard sphere.   Similarly, the link of $x
\in K^e_{\mathcal T}$ is a standard round circle, which is thus large.
The piecewise Euclidean metrics are non-singular at points where $k=2,3$.

When $k=1$, recall that metric singularities of $M^3_{\mathcal T}$ are concentrated on the 1-skeleton 
of $K^e_{\mathcal T}$: there are nonsingular vertical edges in each Euclidean  upper half space,
with trivially large links. Such edges do not lie in $K^e_{\mathcal T}$.
For edges of polygons $P^e_{i,j}$,  the link 
is a 2-sphere with two antipodal distinguished points corresponding to the directions along the edge, connected by   $d$
 spherical geodesic arcs of length $\pi$, where $d$ is the degree of the edge in 
 $K^e_{\mathcal T}$. These arcs divide the sphere into $d$ 2-gons, each having the spherical geometry of a hemisphere. The CAT(1) condition is trivially satisfied, since $d\geq 3$.
 Considered as a point in
 $K^e_{\mathcal T}$, $x$  has link which is a discrete set of $d$ points, and so trivially large.

The potentially non-CAT(0) links are for $x$    a 0-cell.  
In the following, we assume $n=3$, but use $n$ to indicate how our construction yields CAT(0) metrics in higher dimension.
The link  of $x$ in
 $M^n_{\mathcal T}$ is a union of ($n-1$)-dimensional spherical hemispheres $Hem_{x,y}^{n-1}$, one for each horotorus $H_y$ incident at $x$.
 The equatorial  sphere (circle) of each hemisphere is a unit sphere $S^{n-2}_{x,y}$, and is a union of 
 spherical polyhedra $S^{n-2}_{x,y,j}$, each the link of the vertex $x$ in a Euclidean polyhedron $P^e_{y,j}$.
 These are   circular arcs when $n=3$, with 
  length   equal to the angle at a vertex of   $P^e_{y,j}$ incident at $x$, and which 
 add to $2\pi$, giving $S^1$ metrically. Since each  $P^e_{y,j}$ is uniquely identified with another
  $P^e_{y',j'}$,
  the link of $x\in K^e_{\mathcal T}$ is  obtained by identifying the  spheres $S^{n-2}_{x,y}$ along corresponding
   spherical polyhedra $S^{n-2}_{x,y,j}$, $S^{n-2}_{x,y',j'}$ (arcs when $n=3$, giving  a graph).
   
   Rivin shows that the polar dual of an ideal convex hyperbolic polyhedron in 3-space admits a piecewise spherical geometry obtained by gluing together spherical hemispheres along arcs in their boundary circles. This gives a topological 2-sphere containing an embedded graph whose complementary regions are metric hemispheres. All geodesic loops on the graph have lengths at least 2$\pi$: more formally \cite{HR},
for each convex ideal polyhedron $X$ in $\H^3$, let $X^*$ denote the the Poincar\'e dual
of $X$. Assign to each edge $e^*$ of $X^*$ the weight $w(e^*)$ equal to the exterior dihedral
angle at the corresponding edge $e$ of $X$.  
\begin{thm}
 {\rm (Rivin \cite{Ri2})}. The dual polyhedron $X^*$ of a convex ideal polyhedron $X$ in
$\H^3$ satisfies the following conditions:
\begin{description}

\item[ Condition 1.]  $0 < w(e^*) < \pi $  for all edges $e^*$ of $X^*$.
\item[ Condition 2.]  If the edges $e^*_i ,e^*_2, \dots, e^*_k$ form the boundary of a face of $X^*$,
then $w(e^*_1)+w(e^*_2)+  \cdots + w(e^*_k) = 2\pi$.
\item[ Condition 3.]  If $e^*_i ,e^*_2, \dots, e^*_k$ form a simple circuit which does not bound
a face of $X^*$, then $w(e^*_1)+w(e^*_2)+  \cdots + w(e^*_k) > 2\pi$.
 
\end{description}
\end{thm}

This result suffices to prove that $K^e_{\mathcal T}$, and hence $M^3_{\mathcal T}$, is 
metrically CAT(0).  It remains to reconcile the notions of Poincar\'e duality and polar duality.
We recall the following from \cite{HR, Ri1, Ri2, CD} in the notation of the  latter:
Let $X$  denote a convex polyhedron in $\H^3$, viewed in the  hyperboloid model in Minkowski space. 

\defn The \emph{polar
dual} $P(X)$ for $X$ is the set of outward-pointing unit normal vectors to the supporting
hyperplanes of $X$.

\medskip

For a compact polyhedron, each vertex $v$ of $X$ contributes a spherical 
 polyhedron $lk(v)^*$ to $P(X)$:
the
intrinsic metric  on  the polar dual  $P(X)$ is obtained by gluing together the spherical
polyhedra $lk(v)^*$  dual to the vertices $v$ of $X$, by isometrics of their edges in the
combinatorial pattern described above. For ideal polyhedra, such spherical polyhedra are missing
from $P(X)$, which is now a piecewise spherical $(n-2)$-complex with distinguished cycles corresponding to the 
 boundaries of missing $(n-1)$-cells. These are the analogues of edges in Condition 2 of Rivin's 
 characterization above.
    
\medskip 

$P(X)$   inherits the structure of a piecewise spherical cell complex 
as a subset of the de Sitter sphere
in Minkowski space.
Each face  of $X$ contributes a spherical cell for  $P(X)$.
 Rivin proved in his thesis  that the polar dual of an ideal  convex polytope in
hyperbolic 3-space is large: 
Charney and Davis \cite{CD} proved the analogous result in higher dimension. The results of \cite{CD} are more general than what is required here: our interest is in a strict generalization of Rivin's, where $X$
arises as a finite volume ideal polyhedron in the Epstein--Penner construction. 
The pertinent results of \cite{CD} are  Theorem 4.1.1 and Corollary 
4.2.3, which we combine as:

\begin{thm}{\rm (Charney--Davis).}  
Suppose $X$ is a hyperbolic polyhedral set of dimension $n$. Then:
\begin{enumerate}
\item its polar dual $P(X)$ is large;
\item if $\gamma$ is any closed local geodesic of length $2\pi$, then $\gamma$ must lie in the 
subcomplex $P_\gamma$ for some cusp point $y$ of $X$.
\item its completed polar dual $\hat  P(X)$ is large.
\end{enumerate}
\end{thm}

We will explain this notation shortly: that both $M^3_{\mathcal T}$ and $K^e_{\mathcal T}$ are CAT(0) is  then  a consequence  of:

\begin{thm}
The link of the 0-cell $x$ in $K^e_{\mathcal T}$ is the  polar dual $P(X)$
 of the corresponding
dual ideal polyhedron $X$ in the
Epstein--Penner canonical decomposition.
The link of a 0-cell $x$ in $M^3_{\mathcal T}$ is the completed polar dual $\hat P(X)$\end{thm}

\proof
Consider the set $Y$ of parabolic fixed points whose horospheres $\{ H_y\} \subset {\mathcal H}$ meet to define the 0-cell $x$. Then $Y$ is the set of ideal points (cusp points) for the Epstein--Penner canonical ideal polyhedron $X$ dual to the 0-cell $x\in  K_{\mathcal H}$.

\defn For $y\in Y$, let $E_y$  denote the intersection of $X$ with a small horosphere at $y$, and let $P_y$ denote its Euclidean polar dual.

\medskip

Charney--Davis prove that $P_y$ is locally convex in $P(X)$. In our case, $E_y = P^e_{i,j}$ for some ${i,j}$, since the parabolic fixed point $y$ corresponds to $\infty_i$ for some $i$.
The polar dual of a compact  convex Euclidean polyhedron is geometrically  a unit sphere, 
subdivided into
spherical sub-polyhedra. In dimension 2, our situation, 
the polar dual of a compact  convex Euclidean polygon is geometrically the unit circle, 
subdivided into
arcs corresponding to the vertices of the polygon, which measure the external `turning angles'. In the upper half space model $\U\H^3_i$, the faces of the ideal polyhedron $X$ meeting $\infty_i$
are vertical, and their normal directions are horizontal and normal to the Euclidean edges of the polygon $P^e_{i,j}$
obtained by intersecting $X$ with a horosphere $HoP_a$ for $a= 1$, which is  a small enough horosphere.  The collection of polar duals $P_y$ assemble to create $P(X)$:
 
 \begin{lem}
{\rm (\cite{CD}, 2.5.2)} If $y$ is a cusp point,  
 the subcomplex $P_y$ of $P(X)$ corresponding to a cusp is isometric to
the polar dual of a convex set $E_y$  in $E^{n -1}$. 
\end{lem}
 
All vertices of $X$ in the Epstein-Penner construction are ideal vertices, and so all faces of $X$,
and hence all normals to faces of $X$, feature in some subcomplex  $P_y$. Thus $P(X)$ is obtained from the disjoint union $\cup_{y\in Y} \, P_y$. In the upper half space model,   $E_y$ is
the horizontal slice through a cusp end of $X$. These Euclidean polyhedra produce a  tessellation of $HoP_1$, which is geometrically dual to the tessellation by $P^e_{i,j}$. Thus the 
Poincar\'e dual to the polyhedron $X$ yields a Poincar\'e dual to the polyhedron $E_y$ as the `boundary' of $lk(y)$: the corresponding Poincar\'e dual cell decomposition of $\partial E_y$ is combinatorially identical to the link of $x$ in $K^e_{\mathcal H}\cap HoP_1$. Geometrically, each vertex $v$ of $E_y$ contributes a spherical cell to $P_y$, and a spherical cell to the link of $x$, and these are identical. 
These spherical  cells of $P_y$  -- circular arcs in
the case of $P_y$ a circle --  are identified pairwise
in constructing both the link of $x$ in $K^e_{\mathcal H}\cap HoP_1$, and $P(X)$. 
It is now a simple matter to complete the picture of the link of $x$ in  $M^3_{\mathcal T}$:

\defn {\rm (\cite{CD},  4.2.2)}. Let Cone($P_y$) denote the orthogonal join of $P_y$ with a
point. The \emph{completed polar dual} of $X$, denoted $\hat P(X)$, is the piecewise spherical
complex formed by gluing Cone($P_y$) to $P(X)$ along $P_y$ for each cusp
point $y$.

\medskip

In our case,   $X$  has finite volume, with each $P_y$   isometric
to the round sphere $S^{n-2}$,    and so Cone($P_y$) is a
hemisphere. Thus, $\hat P(X)$ is obtained from $P(X)$ by `capping off'  each $P_y$
with hemispheres. It is then clear that $\hat P(X)$ is homeomorphic to the ($n -1$)-sphere,
and moreover geometrically gives the link of $x$ in  $M^3_{\mathcal T}$. Corollaries 4.3.1 and 4.2.3
of \cite{CD} respectively assert that  $ P(X)$ is large, and  $\hat  P(X)$ is large. Accordingly, the techniques of this paper, and the results of Epstein--Penner, Rivin, and Charney--Davis, show that all non-compact finite volume hyperbolic $n$-manifolds admit CAT(0) structures with universal cover obtained as a union of half-spaces, and have CAT(0) spines.

\section{Canonical deformation from hyperbolic to CAT(0) metrics}\label{sec:defo}

In this section we describe how hyperbolic polyhedra admit canonical
metric deformations through hyperbolic polyhedra of constant curvature, limiting on a Euclidean structure. We thank Norman Wildberger for a helpful remark  on  an earlier version of this section, and accordingly define:

 \begin{defn} 
The {\it Wildberger transformation} 
$W_{\tau } \, : \U\H^{n}  \longrightarrow \U\H^{n}, \ \   \tau \geq 0 $,
 is defined by the formula
$$W_{\tau} (({\bf x}, x_{n}))\ :=\ 
({\bf x}, \sqrt{x^2_{n} +\tau^2 }).$$
\end{defn}
Thus $W_{0   }$ is the identity, and
$W_{ \tau }(\U\H^n  ) = \{ ({  x}, x_{n}) \, | \, x_n>\tau\},$
 the region above the horosphere $HoP_{\tau}$.
The transformation  $W_{\tau }$   sends  horizontal horospheres to 
  horizontal horosphere,  and commutes with vertical projection from 
$ \infty$. Wildberger transformations are not isometries, but they preserve the collection of piecewise geodesic subsets:

\medskip

\begin{prop}
If ${\mathcal  C}$ is a 
geodesic arc
in $\U\H^{n} $,
so is  $
W_{\tau}({\mathcal  C}) $.
The hyperplane  $HyS_{y,h}$ is mapped by 
$W_{\tau}$ injectively to  an open disc in the hyperplane   
$HyS_{y,h}$.

\end{prop}

\begin{figure}[h] 
   \centering
   \includegraphics[width=4in]{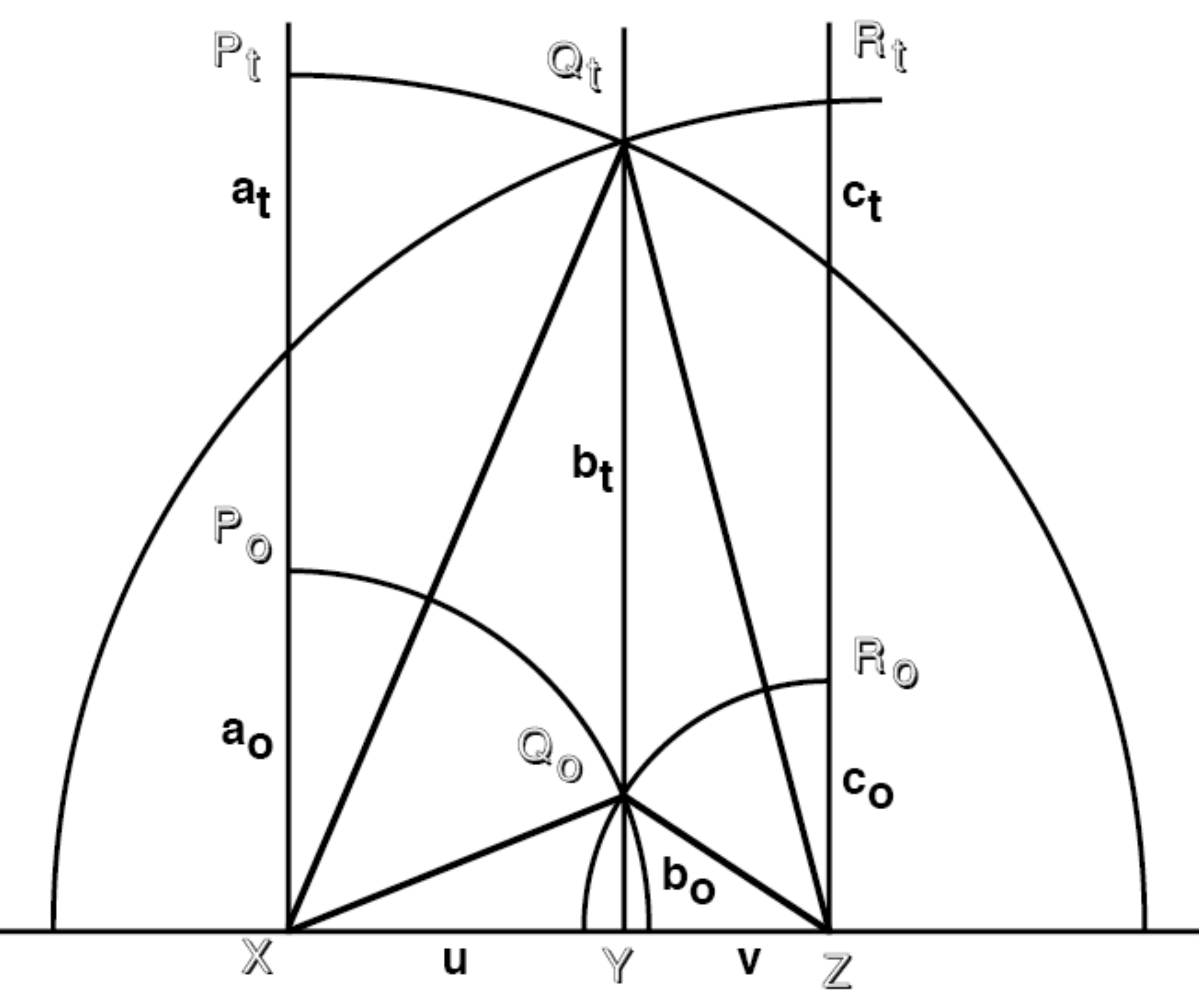} 
   \caption{Wildberger}
   \label{fig:Wildberger}
\end{figure}

\proof This is  a simple calculation using Pythagorus' theorem: 
it suffices to prove this for the hyperbolic plane (the case
$n=2$), since any geodesic arc with endpoints in $\R^{n-1}_\infty$ lies in a vertical plane. This is   illustrated in Figure \ref{fig:Wildberger}, with two given geodesic arcs
$P_0Q_0,\ Q_0R_0$.
Suppose $P_0=(X,a_0),\, Q_0=(Y,b_0),\, R_0=(Y,c_0) \in \U\H^{2}$ are three arbitrary points, with $P_0Q_0$ and $Q_0R_0$ hyperbolic 
geodesic arcs, manifest as arcs of semicircles centred at points on $\R^1_0 := \{({\bf x},0)\} \subset \R^2$. Let   $|XY| = u,\ |YZ|=v$. With respect to the  Euclidean metric on $\R^2$, $|XQ_0| = a_0, \ |Q_0Z|  = c_0$. Thus
$$u^2  + b_0^2 = a_0^2  , \qquad v^2  + b_0^2 = c_0^2$$
\begin{equation}\label{equation:eqPyth}
\Longrightarrow \qquad u^2  + (b_0^2+\tau^2) = (a_0^2+\tau^2)  , \qquad v^2  + (b_0^2+\tau^2) = (c_0^2+\tau^2)
\qquad \forall \ t\in \R \qquad  \qquad 
\end{equation} 

Setting 
$P_\tau = W_{ \tau}(P_0)= (X,a_\tau) ,
\ Q_\tau = W_{ \tau}(Q_0)=(Y,b_\tau),
\ R_\tau = W_{ \tau}(R_0)= (Z,c_\tau),$
we see that Equation \ref{equation:eqPyth} is  equivalent to
$$u^2  + b_\tau^2 = a_\tau^2  , \qquad v^2  + b_\tau^2 = c_\tau^2,$$
and so all points of the hyperbolic geodesic arcs $P_0Q_0,\ Q_0R_0$ move vertically
to corresponding points on the hyperbolic geodesic arcs $P_\tau Q_\tau,\ Q_\tau R_\tau$, which also appear as
arcs of  Euclidean circles: $W_{   \tau}(P_0Q_0) = P_\tau Q_\tau$ as a 
transformation of hyperbolic geodesic arcs, 
as claimed. Every point on an arc moves vertically under $W_{ \tau}$ towards $\infty$,
by a distance depending only on its initial height.

\bigskip

\begin{cor} The image $W_{ \tau}(C)(\Pi)$ 
of any $k_i$-dimensional {hyperbolic} hyperplane $\Pi$ is an open ball in some
$k_i$-dimensional {hyperbolic} hyperplane $\Pi'$.
If $C\subset  \U\H^n $  is contained in some
$k_i$-dimensional {hyperbolic} hyperplane, so is
  $W_{\tau}(C)$.
\end{cor}

\medskip

We now examine how the hyperbolic geometry is distorted.
Consider the angles $\angle P_0Q_0R_0$,
 $\angle P_\tau Q_\tau R_\tau $ between the two arcs at $Q_0,\, Q_\tau$. 
 From basic geometry we have:
   
 \begin{lem} 
 \begin{itemize}
\item  $\angle P_0Q_0R_0\ >\  \angle P_\tau Q_\tau R_\tau $;
\item $\lim_{\tau \to \infty}\,  \angle P_\tau Q_\tau R_\tau  = \pi.$
\item The hyperbolic lengths $|P_\tau Q_\tau |$ satisfy $\lim_{\tau \to \infty}\,  |P_\tau Q_\tau | = 0 $.
\end{itemize}
\end{lem}
 
\begin{figure}[h] 
   \centering
   \includegraphics[width=4in]{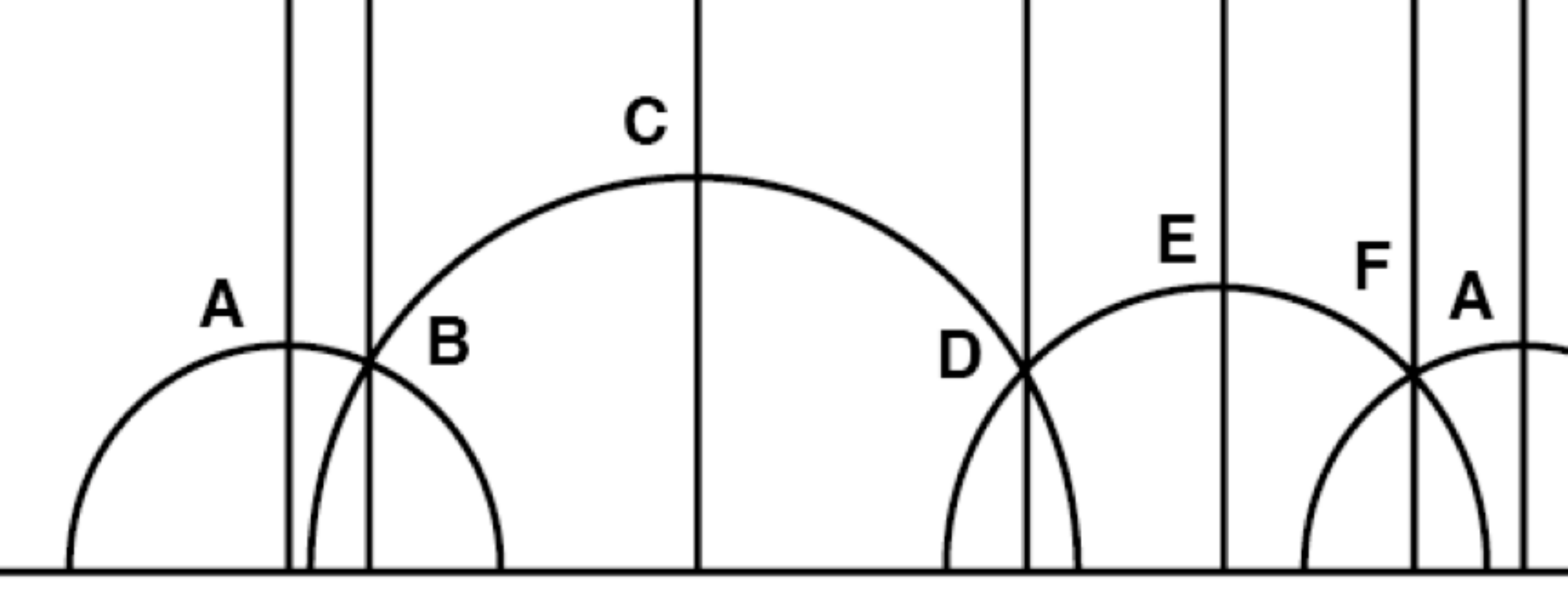} 
   \caption{In the limit, Wildberger transformations applied to a closed neighbourhood of a cusp of a Rieman surface, with piecewise geodesic boundary, give a Euclidean product structure isometric to $S^1\times [1,\infty)$.}
   \label{fig:periodicity2}
\end{figure}

 In Figure \ref{fig:periodicity2} we depict several geodesic arcs $AB,BC,CD,DE,EF,FA$ in the hyperbolic plane with points labeled $A$
 on the left and right to be identified by horizontal translation: the region above these arcs then becomes a neighbourhood of a cusp point on a Riemann surface. Under $W_{\tau}$, these arcs shrink towards the cusp,
 and in the limit approximate a horocycle circle arbitrarily closely; similarly the region above the arcs approximates  
 $S^1\times [\tau,\infty)$ with Euclidean geometry arbitrarily closely as it shrinks and disappears in the limit.
 Note that dilation of the upper plane centred at a point on $\R^1_\infty$ is a hyperbolic isometry, and so we can rescale the picture simultaneously so that $C$ maintains the same height: doing so, the remaining arcs limit to arcs in the horosphere containing $C$ as they shrink to $C$. 
 
Another way to  see this limiting Euclidean geometry, up to scale, is to observe that 
for any (non-vertical) hyperbolic arc $\alpha$, the vertical projection of $W_\tau(\alpha)$ to the horosphere $\R^1_\infty \times \{1\}$ remains constant, giving a Euclidean arc $\alpha_\E$.
  Vertical translation of the Euclidean geometry on the region above $W_\tau(\alpha)$ is an isometry, and so we may take as limiting geometry the union of the metric products $\alpha_\E \times [1,\infty)$, with $\alpha$
  any of the arcs of Figure \ref{fig:periodicity2}.
  
Suppose $M^n$ is any complete hyperbolic $n$-manifold, of finite volume with $p\geq 1$ cusps.
\begin{thm}
For any complete set ${\mathcal T}$ of cusp tori,  transformations $W_\tau$ define a canonical 1-parameter deformation of the hyperbolic geometry
of $M^n$, with limit the canonical piecewise Euclidean structure $M^n_{\mathcal T}$. For fixed $\tau$,
the metric   $M^3_{\tau,{\mathcal T}}$ is a complete, singular piecewise hyperbolic metric with non-singular cusps.  All metric singularities are  concentrated on the $(n-2)$-skeleton of a
piecewise hyperbolic spine 
$K^\tau_{\mathcal T}$.

\end{thm}

\proof $M^n$ is obtained as a quotient  from copies of Ford balls $ FB_{{\mathcal H},i}$ by pairwise identification of hyperbolic 
polyhedra in various $ K^{\infty_i}_{\mathcal H}$. 
Similarly, $M^n{\mathcal T}$ is obtained as a quotient from Euclidean Ford balls $FB^e_{{\mathcal H},i}$ by pairwise identification of Euclidean 
polyhedra in $ K^e_{\mathcal H}$. 

In each $\U\H_i$, we simultaneously apply $W_\tau$. If $A\subset K^{\infty_i}_{\mathcal H}$ is identified by isometry with 
$B\subset K^{\infty_j}_{\mathcal H}$,
then $W_\tau (A) $ and $W_\tau (B) $ continue to be isometric polyhedra of curvature ${-1}$,
although smaller in size. Thus all combinatorial identifications continue to be geometrically feasible by isometry, and equivariantly with respect to group actions: this defines a new metric on $M^n$ for each $\tau, {\mathcal T}$.

By rescaling the curvature by $\tau$ we obtain a  metric deformation through piecewise constant curvature metrics with limit the Euclidean metric $ M^n_{\mathcal T}$.

\section{Weighted points: reconstitution of geometric structure}\label{sec:recon}

 Given a pair $M^3, \ {\mathcal T}$, we have constructed a polyhedral decomposition, with polyhedra in 
 $K^h_{\mathcal T}$ assigned the label identifying hyperplanes in which they lie in $\U\H^3$. 
 These hyperplanes arise as the intersection of a descending horosphere $HoP_{e^-t}$ and an  expanding horosphere $HoS_{*,he^t}$.
 All polygons  $P^e_{i,j} $ arise from the  collision locus of expanding Euclidean circles, as viewed from infinity projected to $HoP_1$.
 Thus knowing the initial moment of birth of each  circle, and its location, the geometric and combinatorial data can be reconstructed for both the  hyperbolic or Euclidean structure: the Ford balls can be reconstructed, as can the combinatorial structure of their boundaries. In this section, we merely record the nature of expansion of circles corresponding to a labeling of polygons:

 \begin{lem}
Suppose the  hyperbolic plane $HyS_{p,h}$ contains a polygon labeled $(p,h)$.
Let $C_{p,h}(t)$ denote the projection to $HoP_1$ of the intersection $S^1_{pq}(t)$ of   $HyS_{p,h}$
with $HoP_{e^-t}$. Then the radius $r_{p,h}(t)$ of this circle satisfies
$$
 r^2_{p,h}(t)  = e^{-2t_0}( 1 -    e^{-2(t - t_0)}), \qquad t\geq  t_0.
$$
 \end{lem}
 
 \def\tanh{{\rm tanh}}
\def\cosh{{\rm cosh}}
\def\coth{{\rm coth}}
\def\sinh{{\rm sinh}}
\def\sech{\rm sech}
\proof Again, it suffices to work in $\U\H^2$: we assume $p=0$, and $h = e^{-t_0}$.
The hyperplane is `born' at $t= t_0$ as the plane $HoP_1$ descends at unit speed, starting when $t=0$. 	
Parametrize the semi-circle $HyS_{0,e^{-t_0}}$  by 
$ x= e^{-t_0}\tanh\,  u,\ y = e^{-t_0}\sech\,  u =  e^{-t}.$
Simple algebra gives the result.

\medskip

Hence  the circle expansion slows exponentially quickly.
Given  an arbitrary finite set of weighted points in the plane, and a lattice ${\mathcal L} \cong \Z\oplus \Z$, we can attempt to create a tessellation of the corresponding elliptic curve. 
One circle may be created in the interior of another expanding circle at a later time; moreover, it may fail to expand to meet the larger one, or overtake to create an edge as expansion continues. However, they data given by the creation of $K^e_{\mathcal T}$ ensures a true tessellation by compact polygons occurs, and $K^h_{\mathcal T}$ can be constructed, allowing the hyperbolic geometry of each Ford ball to be realized.

\begin{thm}
Given the weights assigned to polygons $P^e_{i,j}$, we can reconstruct the hyperbolic metric of $M^3$.
\end{thm}

There are additional properties of such weighted-point sets among arbitrary ones, related to Pythagorus' equation.

\section{
Snappea, SnapPy:  very snappy}\label{sec:wysiwyg}

In the cusped case, the canonical deformations described above
can be applied simultaneously to all fundamental regions of hyperbolic space
itself: the universal cover of a flattened manifold is a flattening of
hyperbolic space, and hence offers a model for hyperbolic space. We see a union
of half spaces, each with boundary plane biperiodically decomposed as a union
of Euclidean polygons. These polygons are pairwise-identified by Euclidean
isometry, and we can therefore imagine navigating in the complement of the
singular 1-skeleton by usual motion in Euclidean space. 
We can consider developing maps into Euclidean space; interesting number
theoretic questions arise concerning the   Euclidean translations and rotations so
obtained.

Both software packages Snappea and SnapPy  \cite{We1, CuDu}
allow the user to see the Ford domains for cusped hyperbolic 3-manifolds, and interactively adjust the defining cusp tori. Thus all of the CAT(0) structure described in this paper is in principle visible in this way.
However, the interface does not provide independent windows for the simultaneous viewing of normalized upper half space models: this would be a valuable addition.

\end{document}